\documentclass[12pt]{article}
\usepackage{amsmath}
\usepackage{amssymb}
\usepackage{latexsym}
\usepackage{array} \newcolumntype{C}{>{$}c<{$}}
\newif\ifscrf\scrftrue
\ifx\footscrfont\nullfont
  \scrffalse
\fi
\newif\iffn\fnfalse

\ifscrf
  \font\footscrfont=rsfs10

\fi

\ifscrf
  \def\Scr#1{\iffn
    \mathchoice{\hbox{\footscrfont #1}}{\hbox{\footscrfont #1}}
    {\hbox{\smallfootscrfont #1}}{\hbox{\tinyfootscrfont #1}}\else
    \mathchoice{\hbox{\scrfont #1}}{\hbox{\scrfont #1}}
    {\hbox{\smallscrfont #1}}{\hbox{\tinyscrfont #1}}\fi}
\else
  \def\Scr{\cal}
\fi

\def\Lop{\Scr L_{\,}}
\def\Lopet{\Scr L_1}
\def\Lopto{\Scr L_2}
\def\Pop{\Scr P_{\,}}

\newcommand{\hovedfont}{\normalfont\bfseries}
\usepackage{theorem}
\theorembodyfont{\normalfont\itshape}
\theoremheaderfont{\hovedfont}
	\theoremstyle{change}
\newtheorem{lemma}{Lemma.}[section]

\newtheorem{prop}[lemma]{Proposition.}
\newtheorem{cor}[lemma]{Corollary.}
	\theorembodyfont{\normalfont}
\newtheorem{eks}[lemma]{Example.}

\newtheorem{BM}[lemma]{Remark.}
\newtheorem{taller}[lemma]{$\!\!$}
	\newenvironment{blanko}[1]%
{\begin{taller}{\hovedfont #1}\normalfont}%
{\end{taller}}

	\newenvironment{bevis}%
{\begin{list}{\em Proof. }%
{\setlength{\labelsep}{0mm}\setlength{\leftmargin}{0mm}%
\setlength{\labelwidth}{0mm}\setlength{\listparindent}{\parindent}%
\setlength{\parsep}{\parskip}\setlength{\partopsep}{0mm}}%
\item}{\qed\end{list}}
	\newenvironment{bevis*}[1]%
{\begin{list}{\em #1 }%
{\setlength{\labelsep}{0mm}\setlength{\leftmargin}{0mm}%
\setlength{\labelwidth}{0mm}\setlength{\listparindent}{\parindent}%
\setlength{\parsep}{\parskip}\setlength{\partopsep}{0mm}}%
\item}{\qed\end{list}}
	\newenvironment{blanko*}[1]%
{\begin{list}{\bf {#1} }%
{\setlength{\labelsep}{0mm}\setlength{\leftmargin}{0mm}%
\setlength{\labelwidth}{0mm}\setlength{\listparindent}{\parindent}%
\setlength{\parsep}{\parskip}\setlength{\partopsep}{0mm}}%
\item}{\end{list}}
\newcounter{dummycounter}
\newenvironment{punkt-i}%
{%
	\begin{list}%
	{(\roman{dummycounter})}%
	{\usecounter{dummycounter}%
	\setlength{\itemsep}{0em}\setlength{\parsep}{0em}\setlength{\topsep}{0em}%
	\setlength{\itemindent}{0em}\setlength{\labelwidth}{1.8em}%
	\setlength{\labelsep}{0.6em}\setlength{\leftmargin}{2.4em}}%
}%
{\end{list}}
\renewcommand{\epsilon}{\varepsilon}
\makeatletter
\renewcommand{\ldots}{\relax\ifmmode\ldotp\ldotp\ldotp%
\else$\m@th\ldotp\ldotp\ldotp\ $\fi}

\newcommand{\comentario}{%
\noindent [GKP] refers to Graber-Kock-Pandharipande~\cite{Grab-Kock-Pand}.
\medskip%
}

\renewenvironment{thebibliography}[1]
     {\section*{\refname
        \@mkboth{\MakeUppercase\refname}{\MakeUppercase\refname}}%
		  \comentario
      \list{\@biblabel{\@arabic\c@enumiv}}%
           {\settowidth\labelwidth{\@biblabel{#1}}%
            \leftmargin\labelwidth
            \advance\leftmargin\labelsep
            \@openbib@code
            \usecounter{enumiv}%
            \let\p@enumiv\@empty
            \renewcommand\theenumiv{\@arabic\c@enumiv}}%
      \sloppy
      \clubpenalty4000
      \@clubpenalty \clubpenalty
      \widowpenalty4000%
      \sfcode`\.\@m}
     {\def\@noitemerr
       {\@latex@warning{Empty `thebibliography' environment}}%
      \endlist}

\makeatother
\newcommand{\brac}[1]{\,\langle \, {#1} \, \rangle\,}
\newcommand{\bbrac}[1]%
{\,\langle\!\langle \, {#1} \, \rangle\!\rangle\,}
\newcommand{\Brac}[1]{\,\langle \ {#1} \ \rangle\,}
\newcommand{\Bbrac}[1]%
{\,\langle\!\langle \ {#1} \ \rangle\!\rangle\,}
\newcommand{\cp}{ {\scriptstyle \;\cup\;} }

\newcommand{\smallbinom}[2]%
{{\textstyle \binom{#1}{#2}}}
\newcommand{\smallprod}[2]%
{\overset{#2}{\underset{#1}{\textstyle{\prod}}}} 
\newcommand{\smallsum}[2]%
{\overset{#2}{\underset{#1}{\textstyle{\sum}}}} 
\newcommand{\bigprod}[2]%
{\overset{#2}{\underset{#1}{\text{{\Large\(\prod\)}}}}}
\newcommand{\bigsum}[2]%
{\overset{#2}{\underset{#1}{\text{{\Large\(\sum\)}}}}}

\newcommand{\bigfatgreek}[1]{\boldsymbol{#1}} 
\newcommand{\evalclass}{\bigfatgreek\eta} 
\newcommand{\psiclass}{\bigfatgreek\psi} 
\newcommand{\mpsiclass}{\ov\psiclass{}} 
\newcommand{\deltaclass}{\bigfatgreek\delta} 
\newcommand{\boundary}{\bigfatgreek\xi}

\newcommand{\XX}{\mathfrak{X}}
\providecommand{\qed}%
{\hspace*{\fill}\nolinebreak[1]\hspace*{\fill} $\Box$}

\renewcommand{\d}{\partial}

\newcommand{\pil}{\rightarrow}
\newcommand{\langpil}{\longrightarrow}

\newcommand{\isopil}%
{\stackrel{\raisebox{0.1ex}[0ex][0ex]{\(\sim\)}}%
{\raisebox{-0.15ex}[0.28ex]{\(\pil\)}}}

\newcommand{\shortsetminus}%
{\,\raisebox{1pt}{\ensuremath{\mathbb{r}}\,}}
\newcommand{\kan}{\mbox{\Large \(\omega\)}}
\newcommand{\OO}{\ensuremath{{\mathcal O}}}

\newcommand{\upperstar}{^{\raisebox{-0.25ex}[0ex][0ex]{\(\ast\)}}}
\newcommand{\lowerstar}{_{\raisebox{-0.33ex}[-0.5ex][0ex]{\(\ast\)}}}
\newcommand{\df}%
{\: {\raisebox{0.255ex}{\normalfont\scriptsize :\!\!}}=}
\newcommand{\tensor}{\otimes}
\renewcommand{\dim}{\operatorname{dim}}
\newcommand{\codim}{\operatorname{codim}}

\newcommand{\brok}[2]{{\textstyle{\frac{#1}{#2}}}}
\newcommand{\ov}{\overline}

\newcommand{\fat}[1]{\mathbf{{#1}}}
\renewcommand{\P}	{\mathbb{P}}
\newcommand{\A}		{\mathbb{A}}

\newcommand{\Z}		{\mathbb{Z}}
\newcommand{\Q}		{\mathbb{Q}}

\newcommand{\mtau}{\ov\tau}
\newcommand{\mbtau}{\overline{\boldsymbol{\tau}}}

\newcommand{\cusp}{\mathrm{K}}
\newcommand{\cuspl}{\mathrm{Kl}}
	\newcommand{\CL}{C\!L}
	\newcommand{\KL}{K\!L}
\newcommand{\cuspp}{\mathrm{Kp}}
	\newcommand{\CP}{C\!P}
	\newcommand{\KP}{K\!P}
\newcommand{\flex}{\mathrm{F}}

\newcommand{\lineartwiginH}{\mathrm{I}}
	\newcommand{\NF}{N\!F}

\renewcommand{\boundary}{\bigfatgreek\xi}

\begin{document}
  

\vspace*{0.5in}
\begin{center}
   {\LARGE {Characteristic numbers of rational curves \\
	with cusp or prescribed triple contact} \par}%
   \vskip 1.5em
   {\large \lineskip .5em
      Joachim Kock\footnote{%
      Supported by the National Science Research Council of Denmark.

      E-mail address: {\tt jojo@dmat.ufpe.br}}
   \par}%
   \vskip 1em
   {\normalsize Universidade Federal de Pernambuco\\
   Recife, Brazil}%
\end{center}%
\par
\vskip 1.5em

\begin{abstract}
  This note pursues the techniques of modified psi classes on the stack of
  stable maps (cf.~[Graber-Kock-Pandharipande]) to give concise solutions to the
  characteristic number problem of rational curves in $\P^2$ or
  $\P^1\times\P^1$ with a cusp or a prescribed triple contact.  The classes of
  such loci are computed in terms of modified psi classes, diagonal classes, and
  certain codimension-2 boundary classes.  Via topological recursions the
  generating functions for the numbers can then be expressed in terms of the
  usual characteristic number potentials.
\end{abstract}

\section*{Introduction}

With the advent of stable maps and quantum cohomology
(Kontsevich-Manin~\cite{KM:9402}), there has been a tremendous progress in
enumerative geometry.  One subject of much research activity has been the
characteristic number problem, notably for rational curves. 
Highlights of these developments include Pandharipande~\cite{Pand:9504}, who
first determined the simple characteristic numbers of rational curves in
projective space; Ernstr\"om-Kennedy~\cite{EK1} who computed the numbers for
$\P^2$ using stable lifts --- a technique that also allowed to determine
characteristic numbers including a flag condition, as well as characteristic
numbers of cuspidal plane curves; and Vakil~\cite{Vakil:9803} who used
degeneration techniques to give concise recursions for the characteristic
numbers also for elliptic curves.

With the notions of modified psi classes and the tangency quantum potential
introduced in Graber-Kock-Pandharipande~\cite{Grab-Kock-Pand}, conceptually
simpler solutions were given to the characteristic number problem for rational
curves in any projective homogeneous space, as well as for elliptic curves in
$\P^2$ or $\P^1\times\P^1$.  Tangency conditions allow simple expressions in
terms of modified psi classes, and then the solutions follow from standard
principles in Gromov-Witten theory, e.g.~topological recursion.

Having settled the question of characteristic numbers of nodal rational curves,
a natural next problem to consider is that of cuspidal curves, or to impose
higher order contacts, e.g.~specified flex lines.  Schubert~\cite{Schubert}
computed the characteristic numbers of cuspidal plane cubics, and in the
1980's, a lot of work was devoted to the verification of his results,
cf.~Sacchiero~\cite{Sacchiero},
Kleiman-Speiser~\cite{Kleiman-Speiser:cuspidal-cubics},
Miret-Xamb\'o~\cite{Miret-Xambo:complete-cuspidal-cubics}, and
Aluffi~\cite{Aluffi:nodal-cuspidal}.

The techniques of stable lifts allowed L.~Ernstr\"om and G.~Kennedy~\cite{EK1}
to determine the characteristic numbers of plane rational curves with cusp for
{\em any} degree, and their joint paper with S.J.~Colley~\cite{CEK} represents a
big advance in the treatment of third order contacts.

\bigskip

The present note shows how the techniques of modified psi classes developed in
Graber-Kock-Pandharipande~\cite{Grab-Kock-Pand} (henceforth cited as [GKP]), can
also be used to solve the characteristic number problem for cuspidal rational
curves in $\P^2$ or $\P^1\times\P^1$ (as well as that of a single triple
contact).  To this end, a slight generalisation of the first enumerative
descendants is needed, namely allowing a single factor of the top product to be
a square of a modified psi class or a certain codimension-2 boundary class.  Via
topological recursions, the corresponding potentials are related to the usual
tangency quantum potential.  The locus of cuspidal curves and the locus of
curves with triple contact are described in terms of these classes, whereafter
the differential equations translate into equations for the sought
characteristic numbers.

The constructions and computations do not pretend to be particularly deep or
difficult; the raison d'\^etre of this note is rather to illustrate the
versatility of the methods developed in [GKP].  
Results and notation from that paper are briefly reviewed in 
Section~\ref{Sec:prelim}, and in \ref{setupG} and \ref{P1xP1}.


The material of this note constituted Chapter~4 of my PhD
thesis~\cite{Kock:thesis}, and it is a pleasure here to thank the Departamento
de Matem\'atica da Universidade Federal de Pernambuco for four lovely years, and
in particular my advisor Israel Vainsencher for his guidance and encouragement. 
I have also profited from conversations with Letterio Gatto and Lars Ernstr\"om.

\section{Preliminaries}
\label{Sec:prelim}

\begin{blanko}{The target space.}  
Throughout we work over the field of complex numbers.  Let $X$ denote a
projective homogeneous variety, and let $T_0,\ldots,T_r$ denote the elements of
a homogeneous basis of the cohomology space $H\upperstar (X,\Q)$.  In the
applications $X$ will be $\P^2$ or $\P^1\times\P^1$.  Let $g_{ij}$ denote the
Poincar\'e metric constants $\int_X T_i\cp T_j$;  we set also $g_{ijk} = \int_X
T_i \cp T_j \cp T_k$.  Let $(g^{ij})$ be the inverse matrix to $(g_{ij})$.  It
is used to raise indices as needed; in particular, with $g_{ij}^k = \sum_e
g_{ije}\,g^{ek}$, we can write $T_i \cp T_j = \smallsum{k}{} g_{ij}^k \; T_k$.
\end{blanko}

\begin{blanko}
	{The deformed metric.} (See Kock~\cite{Kock:0006} for details.)  Let $\fat y
	= (y_0,\ldots,y_r)$ be formal variables, and put 
	$$
	\phi = \sum_{\fat s} \frac{\fat
	y^{\fat s}}{\fat s!} \int_X \fat T^{\fat s} ,$$ 
	with usual multi-index notation, $\fat s!=s_0!\cdots s_r!$, $\fat y^{\fat s}
	= y_0^{s_0}\cdots y_r^{s_r}$, and $\fat T^{\fat s} = T_0^{s_0}\cp\cdots\cp 
	T_r^{s_r}$.	
	Consider its partial derivatives $\phi_{ij} =
	\sum_{\fat s} \frac{\fat y^{\fat s}}{\fat s!} \int_X \fat T^{\fat s}\cp T_i
	\cp T_j$, and use the matrix $(g^{ef})$ to raise indices, putting
\begin{equation}\label{gandphi}
 \phi^i_j = \smallsum{e}{} g^{ie} \, \phi_{ej}  ,
\quad\text{ and }\quad 
 \phi^{ij} =\smallsum{e,f}{} g^{ie} \, \phi_{ef} \, g^{fj}.
\end{equation}
The entities $\phi^i_j(\fat y)$ 
are the tensor elements of `multiplication by the exponential', precisely
\begin{align}
\smallsum{\fat s}{}\frac{\fat y^{\fat s}}{\fat s!} \fat T^{\fat s} \cp T_p 
\label{exptophi}
 & = \smallsum{e}{} \, T_e \, \phi^e_p(\fat y) .
\end{align}
The {\em deformed metric} is the non-degenerate
symmetric bilinear pairing $H\upperstar 
	(X,\Q) \pil \Q[[\fat y]]$ given by the tensor elements
	$$
	\gamma_{ij} \df 
	\phi_{ij}(-2\fat y).
	$$
	The inverse of the matrix $(\gamma_{ij})$ is given by
	$$
	\gamma^{ij}
	= \phi^{ij}(2\fat y)  = \smallsum{e,f}{} \phi^i_e g^{ef} \phi_f^j .
	$$
	We will also need certain derivatives of this,
	$$
	\gamma_k^{ij} \df \phi_k^{ij}(2\fat y) = \smallsum{e,f}{} \phi^i_e g_k^{ef}
	\phi_f^j .
	$$
\end{blanko}

\begin{blanko}{Moduli of stable maps.}
  Let $\ov M_{0,S}(X,\beta)$ denote the moduli stack of Kontsevich stable maps
  of genus zero whose direct image in $X$ is of class $\beta\in H_2^+(X,\Z)$,
  and whose marking set is $S=\{p_1,\ldots,p_n\}$.  For each mark $p_i$, let
  $\nu_i : \ov M_{0,S}(X,\beta) \pil X$ denote the evaluation morphism.  The
  reader is referred to Fulton-Pandharipande~\cite{FP-notes} for definitions and
  basic properties of stable maps, Gromov-Witten invariants and quantum
  cohomology.
\end{blanko}

\begin{blanko}{Modified psi classes and diagonal classes} 
(cf.~[GKP].)  The psi class $\psiclass_i$ is the first Chern class of the line
bundle on $\ov M_{0,S}(X,\beta)$ whose fibre at a moduli point $[\mu:C\to X]$ is
the cotangent line of $C$ at $p_i$.  On a moduli space $\ov
M_{0,S}(X,\beta)$ with $\beta>0$, let $\boundary_i$ denote the sum of all
boundary divisors classes having mark $p_i$ on a contracting twig.  The {\em
modified psi class} is defined as
   $$
   \mpsiclass_i \df \psiclass_i - \boundary_i.
   $$
A crucial observation is that $\mpsiclass_i$ is invariant under 
pull-back along forgetful morphisms.

The $ij$'th diagonal class $\deltaclass_{ij}$ is by definition the sum of all
boundary divisor classes having $p_i$ and $p_j$ together on a contracting twig.
%
The diagonal classes appear as correction terms when restricting a modified psi
class to a boundary divisor $D$ both of whose twigs are of positive degree.  
If $D$ is the image of the gluing morphism
$$
\rho_D : \ov M_{0,S' \cup \{x\}}(X,\beta') 
\times_X \ov M_{0,S'' \cup \{x\}}(X,\beta'')
\langpil \ov M_{0,S}(X,\beta)
$$
then 
\vspace*{-1em}
\begin{equation*} \label{restr-mpsi}
\rho_D\upperstar \ov \psiclass_i = 
\mpsiclass_i+\deltaclass_{ix},
\end{equation*} 
where $x$ denotes the gluing mark.
	\end{blanko}

%
%
	
\begin{blanko}{The tangency quantum potential} \label{Gamma}
(cf.~[GKP]).
	The integrals 
   $$
   \Brac{\ov\tau_{k_1}(\gamma_1)\cdots \ov\tau_{k_n}(\gamma_n)}\!_\beta \df
\int \ov\psiclass{}_1^{k_1}\cp\nu_1\upperstar (\gamma_1) \cp\cdots \cp
\ov\psiclass{}_n^{k_n}\cp\nu_n\upperstar (\gamma_n)
\ \cap [\ov M_{0,n}(X,\beta) ]
   $$
	($\gamma_i \in H\upperstar (X,\Q)$) are called {\em enumerative descendants}.
For the {\em first} enumerative descendants (exponent at most $1$ on modified psi classes)
we employ the notation
$$
\brac{\mbtau_0^{\fat{a}} \; \mbtau_1^{\fat{b}}}\!_{\beta} \df 
\Brac{ \smallprod{k=0}{r} 
(\mtau_0(T_k))^{a_k}(\mtau_1(T_k))^{b_k}}\!_{\beta} .
$$
Their generating function
is called the {\em tangency quantum potential}:
$$
\Gamma(\fat{x},\fat{y}) = \sum_{\beta>0}\sum_{\fat{a},\fat{b}} 
\frac{\fat{x}^{\fat{a}}}{\fat{a}!}\frac{\fat{y}^{\fat{b}}}{\fat{b}!}
\brac{\mbtau_0^{\fat{a}} \; \mbtau_1^{\fat{b}}}\!_{\beta} .
$$

	The tangency quantum potential satisfies the topological recursion relations
\begin{equation}\label{PDE}
	\Gamma_{y_k x_i x_j } = \Gamma_{x_k(x_i x_j)} - \Gamma_{(x_k x_i) x_j}
	- \Gamma_{(x_k x_j) x_i} +
	\sum_{e,f} \; \Gamma_{x_k x_e} \; \gamma^{ef} \; \Gamma_{x_f x_i x_j}.
	\end{equation}
	Here, and in the sequel,  subscripts on potentials denote partial derivatives,
	e.g.~$\Gamma_{x_i} \df\frac{\d}{\d x_i} \Gamma$, and we set also
$$
\Gamma_{(x_i x_j)} \df \smallsum{k=0}{r} g_{ij}^k \ \Gamma_{x_k} 
	\qquad \text{ and } \quad
		\Gamma_{(y_i x_j)} \df \smallsum{k=0}{r} g_{ij}^k \ \Gamma_{y_k}
.
$$
\end{blanko}



\section{Slightly enriched first enumerative descendants}
\label{Sec:PQ}

\begin{blanko}{Pi classes.}
Let $\Pi_i$ denote the sum of all codimension-2 boundary classes
whose middle twig has degree 0 and carries the mark $p_i$, while the two other
twigs have positive degree. Clearly $\Pi_i$ is compatible with pull-back along 
forgetful morphisms.
\end{blanko}

\begin{lemma}\label{mpsi.Pi=0}
	Let $\Delta$ be an irreducible component of $\Pi_i$, then
	$$
	 \mpsiclass_i \cdot \Delta= 0.
	$$
\end{lemma}
\begin{bevis}
	Let $x'$ and $x''$ denote the two attachment points on the middle twig of
	$\Delta$.  Now restrict $\mpsiclass_i = \psiclass_i - \boundary_i$ to the
	moduli space $M_i$ corresponding to the middle twig: the psi class
	$\psiclass_i$ restricts to give the corresponding psi class of the mark of
	the middle twig.  Restricting $\boundary_i$ to $M_i$ corresponds to breaking
	off a twig containing $p_i$ but not $x'$ nor $x''$.  In other words, the
	restriction of $\boundary_i$ is $(p_i \mid x' , x'')$ on $M_i$, which is the
	well-known boundary expression for $\psiclass_i$, so altogether the
	restriction of $\mpsiclass_i$ is zero.
\end{bevis}

Let $P^k$ denote the generating function for top products of classes of type
$\mtau_0(T_i)$ and $\mtau_1(T_j)$ and a single factor (say at the first
mark) of type $\Pi_1 \cp \nu_1\upperstar (T_k)$.  Precisely
$$
P^k(\fat{x},\fat{y}) \df \sum_{\beta>0}\sum_{\fat{a},\fat{b}}
\frac{\fat{x}^{\fat{a}}}{\fat{a}!}\frac{\fat{y}^{\fat{b}}}{\fat{b}!} \brac{\Pi_1
\cp \nu_1\upperstar (T_k) \ \mbtau_0^{\fat{a}} \, \mbtau_1^{\fat{b}}}\!_{\beta} .
$$
(In view of Lemma~\ref{mpsi.Pi=0}, there is no reason for allowing also marks 
combining $\Pi$ with $\mpsiclass$.)


\begin{prop}\label{PkGamma}
	The following differential equation relates $P^k$ to $\Gamma$:
	$$
	P^k = \brok{1}{2} \sum_{e,f} \Gamma_{x_e} \; \gamma_k^{ef} \; \Gamma_{x_f}.
	$$
\end{prop}
\begin{bevis}
	Among the components of $\Pi_1$, those with at least two marks on the middle
	contracting twig have zero push-down under forgetting $p_1$, so we need only
	consider components of $\Pi_1$ where $p_1$ is alone on the middle twig.  Each
	such component $\Delta$ is the image of a birational morphism from a triple
	fibred product $\ov M_{0,S\cp \{x'\}}(X,\beta') \times_X \ov M_{0,3}(X,0) 
	\times_X \ov M_{0, S''\cp \{x''\}}(X,\beta'')$.  Therefore there is the
	following sort of splitting lemma, similar to Lemma 1.5 of
	Kock~\cite{Kock:0006}: \setlength{\multlinegap}{0pt}
	\begin{multline*}
		\brac{ \Delta \cdot \nu_1\upperstar (T_k) \cdot
	\mbtau_0^{\fat{a}}\mbtau_1^{\fat{b}}}\!_\beta =
	\underset{p'', q''}{\sum_{p',q'}} \sum_{\fat{s}',\fat{s}''}
	\smallbinom{\fat{b}'}{\fat{s}'}\smallbinom{\fat{b}''}{\fat{s}''}
	\brac{\mbtau_0^{\fat{a}'}\mbtau_1^{\fat{b}'\! - \! \fat{s}'} 
	\mtau_0(\fat{T}^{\fat{s}'}\!\cp T_{p'})}\!_{\beta'}
	\\ g^{p'q'} \; g_{q' k p''} \; g^{p'' q''} \  
	\brac{\mtau_0(\fat{T}^{\fat{s}''}\!\cp T_{q''})
	\mbtau_0^{\fat{a}''}\mbtau_1^{\fat{b}''\! - \! \fat{s}''} }\!_{\beta''}.
	\end{multline*}
	Translating this into a statement about the potentials yields the wanted
	differential equation.  It is perhaps opportune to explain the appearance of
	the splitting factor $\gamma_k^{ef}$.  At the gluing mark $x'$ on the
	one-primed twig there is (after expressing things in terms of potentials) a
	factor $ \sum_{\fat s'}\frac{\fat y^{\fat s'}}{\fat s'}\fat T^{\fat s'} \cp
	T_{p'} = \sum_{e} T_{e} \phi^{e}_{p'}$, cf.~(\ref{exptophi}).  Arguing
	similarly on the two-primed twig we conclude that the splitting factor is
	$$
  \sum \phi_{p'}^e \;
	\big( g^{p'q'} \; g_{q' k p''} \; g^{p'' q''} \big) \; \phi_{q''}^f = \sum \phi_{p'}^e \;
	g^{p'q''}_k \; \phi_{q''}^f = \gamma_k^{ef} .
	$$
	
	Note the presence of the factor $\brok{1}{2}$, due to the fact that all
	the components of $\Pi$ appear twice in the sum, depending on which of the
	two outer twigs we consider to be the one-primed and which is two-primed.  In
	the very special case where $p_1$ is the only mark in play, there is no
	repetition in the sum since nothing distinguishes the two twigs, but this
	very symmetry means that the morphism from the fibred product is actually
	two-to-one, so in this case we divide by two for this reason.
\end{bevis}

Let $Q^k(\fat{x},\fat{y})$ denote the generating function corresponding to first
enumerative descendants allowing a single quadratic modified psi class, say at
the first mark:
$$
Q^k(\fat{x},\fat{y}) \df
\sum_{\beta>0}\sum_{\fat{a},\fat{b}} 
\frac{\fat{x}^{\fat{a}}}{\fat{a}!}\frac{\fat{y}^{\fat{b}}}{\fat{b}!}
\brac{\mbtau_0^{\fat{a}} \, \mbtau_1^{\fat{b}} \ \mtau_2(T_k)}\!_{\beta} .
$$

\begin{prop}\label{QkGamma}
	The following differential equation relates $Q^k$
	to the tangency quantum potential:
	$$
	Q^k_{x_i x_j} = 
	\Gamma_{(x_i x_j)y_k} - \Gamma_{(y_k x_i) x_j} - \Gamma_{(y_k x_j) x_i}
	+ \sum_{e,f} \; \big(\Gamma_{y_k x_e} + \Gamma_{(x_k x_e)}\big) \ \gamma^{ef} 
	\ \Gamma_{x_f x_i x_j}.
	$$
\end{prop}

\begin{bevis}
	The proof is similar to the proof of Equation~(\ref{PDE})
	(see~[GKP], 2.1.1 and
	\cite{Kock:0006}, 3.4.)
	Let the mark $p_1$ correspond to the class $\mtau_2(T_k)$, and let $p_2$ and
	$p_3$ carry the extra classes $\mtau_0(T_i)$ and $\mtau_0(T_j)$ corresponding
	to the derivatives.  Take one of the two modified psi classes $\mpsiclass_1$
	and write it as sum of boundary divisors, to each of which the remaining
	factors are restricted.  The first three terms correspond to boundary
	divisors with trivial degree splitting; compared to Equation~(\ref{PDE}),
	they each have a derivative with respect to $y$ instead of $x$ because there
	is now one modified psi class left on $p_1$.  As to the quadratic term, it
	splits up in two, because the factor $\mpsiclass_1 \cdot \nu_1\upperstar
	(T_k)$ restricts to give $\mpsiclass_1 \cdot \nu_1\upperstar (T_k) +
	\deltaclass_{1x'}\cdot\nu_{x'}\upperstar (T_k)$, sending the evaluation class
	$\nu_1\upperstar (T_k)$ over to the gluing mark $x'$.  This explains the
	factor $\big(\Gamma_{y_k x_e} + \Gamma_{(x_k x_e)}\big)$ in the quadratic
	term.
\end{bevis}

Observe that $\sum \Gamma_{x_e} \gamma_k^{ef} =
\sum \Gamma_{(x_k x_m)} \gamma^{mf}$, so the last quadratic term is very similar 
to the terms of $P^k_{x_i x_j}$.

	For $k=0$, there is a much simpler equation:

\begin{cor}\label{Q0}
	$$
	Q^0 = - \brok{1}{2} \sum_{e,f} \Gamma_{x_e} \, \gamma^{ef} \, \Gamma_{x_f} .
	$$
\end{cor}

\begin{bevis}
	After using the dilaton equation 
	$\Gamma_{y_0} = -2 \Gamma$ twice, the equation of the proposition reads
	$$
	Q^0_{x_i x_i} = 
	-2\Gamma_{(x_i x_i)} - 2 \Gamma_{y_i x_i}
	- \sum_{e,f} \; \Gamma_{x_e} \ \gamma^{ef} \ \Gamma_{x_f x_i x_i}.
	$$
	Now apply topological recursion to the second term
	and simplify, ending up with
	$$
	Q^0_{x_i x_i} = - \sum_{e,f} \; \Gamma_{x_i x_e} \ \gamma^{ef} \ \Gamma_{x_f x_i}
	- \sum_{e,f} \; \Gamma_{x_e} \ \gamma^{ef} \ \Gamma_{x_f x_i x_i}.
	$$
	Integrating twice with respect to $x_i$ yields the result.
\end{bevis}

\begin{BM}\label{P0+Q0}
	It is immediate from the formulae that $P^0 + Q^0 = 0$.
In fact, more generally, the classes $-\Pi_1$ and $\mpsiclass_1^2$ on
one-pointed space $\ov M_{0,1}(X,\beta)$ have the same push-down in $\ov
M_{0,0}(X,\beta)$.  Indeed, generally $\Pi_1$ pushes down to give the 
whole boundary. 
On the other hand, the push-down of $\mpsiclass_1^2 = \psiclass_1^2$ is the
kappa class by definition (see
Arbarello-Cornalba~\cite{Arbarello-Cornalba:coh}), and according to
Pandharipande~\cite{Pand:9504}, Lemma~2.1.2, the kappa class is minus the
boundary.  (That proof treats $\P^r$ but it carries over to the present
case.)
\end{BM}

\section{Cuspidal curves in $\P^2$}
\label{Sec:P2}

  In this section we consider $X=\P^2$, with its usual cohomology basis ($T_0= $
  fundamental class, $T_1 =$ line, $T_2 =$ point).  Set $\evalclass_i \df c_1(
  \nu_i\upperstar (T_1))$.
%
%
\begin{blanko}{The characteristic number potential} (cf.~[GKP] \S 4).\label{setupG}
  Let $N_d(a,b,c)$ denote the number of irreducible plane
  rational curves of degree $d$ which pass through $a$ general points, are
  tangent to $b$ general lines, and are tangent to $c$ general lines at a
  specified point, and define the number to be zero unless $a+b+2c = 3d-1$.
	
	Let $\Omega$, $\Theta$, and $\Xi$ denote the classes corresponding to these
	three conditions, then (at mark $p_1$, say) we have $\Omega= \evalclass_1^2$,
	$\Theta = \evalclass_1 (\evalclass_1 +\mpsiclass_1) $, and $\Xi =
	\evalclass_1^2 \mpsiclass_1 $.  The
 characteristic number potential
$$
G(s,u,v,w) = \sum_{d>0} \exp(ds) 
\sum_{a,b,c} \frac{u^a}{a!}\frac{v^b}{b!}\frac{w^c}{c!} \; N_d(a,b,c)
$$
is related to the tangency quantum potential $\Gamma$ by 
\begin{align*}
G(s,u,v,w) &= \Gamma(x_1,x_2,y_1,y_2) ,
\end{align*}
subject to the change of variables:
\begin{equation}
\label{cvar}
x_1=s , \qquad  x_2 = u + v ,  \qquad y_1 = v , \qquad y_2 = w .
\end{equation}
For simplicity we set $x_0=y_0=0$ throughout.


For the deformed metric we have
	$$
	(\gamma^{ef}) = \begin{pmatrix}
	\phantom{0}0\phantom{0} & \phantom{0}0\phantom{0}& \phantom{0}1\phantom{0} \\
	0 & 1 & 2y_1 \\
	1 & 2y_1 & 2y_1^2 + 2y_2
	\end{pmatrix} = \begin{pmatrix}
	\phantom{0}0\phantom{0} & \phantom{0}0\phantom{0}& \phantom{0}1\phantom{0} \\
	0 & 1 & 2v \\
	1 & 2v & 2v^2 + 2w
	\end{pmatrix} ,
	$$
so in terms of the two differential operators
	\begin{align*}
\Lop &\;\df\; \phantom{2v}\tfrac{\d}{\d s} + 2v \tfrac{\d}{\d u} ,\\
\Pop &\;\df\; 2v\tfrac{\d}{\d s} + (2v^2+2w) \tfrac{\d}{\d u} ,
\end{align*}	
\end{blanko}
the topological recursion relations satisfied by the characteristic number
potential can be written
\begin{align}
\label{eq1GP2}
G_{v s} &=  G_{u s} - G_{u}
+\brok{1}{2} \big(  G_{s s}\cdot \Lop G_{s}
+ G_{u s}\cdot \Pop G_{s} \big) ,
\\[8pt]
\label{eq2GP2}
G_{w s s} &= G_{u u}
+ \big( G_{u s} \cdot \Lop G_{ss}
+ G_{u u} \cdot \Pop G_{ss} \big) .
\end{align}

\begin{blanko}{The slightly enriched potentials.}\label{PQforP2}
  Combining Propositions~\ref{PkGamma} and \ref{QkGamma} with the above 
  coordinate changes, we can express the slightly enriched potentials
  in terms of the characteristic number potential. We have  
$$
	(\gamma_1^{ef}) = \begin{pmatrix}
	\phantom{0}0\phantom{0} & \phantom{0}0\phantom{0}& \phantom{0}0\phantom{0} \\
	0 & 0 & 1 \\
	0 & 1 & 2y_1
	\end{pmatrix}
	\qquad \text{ and } \qquad
	(\gamma_2^{ef}) = \begin{pmatrix}
	\phantom{0}0\phantom{0} & \phantom{0}0\phantom{0}& \phantom{0}0\phantom{0} \\
	0 & 0 & 0 \\
	0 & 0 & 1
	\end{pmatrix},
	$$
so from Proposition~\ref{PkGamma} we get
\begin{align}
\label{P1forP2}
P^1_{x_1 x_1} 
&= G_{us} G_{ss} + G_{uss} G_s + G_u \!\cdot\! \Lop G_{ss} + G_{us}
\!\cdot\! \Lop G_s,
\\[8pt]
P^2_{x_1 x_1} 
&= G_{us}G_{us} + G_{uss} G_u .
\label{P2forP2}
\end{align}
Here we have taken double derivative with respect to $x_1$, anticipating the
applications.

Similarly, for  the $Q$-potential, Proposition~\ref{QkGamma} gives
these three equations:
\begin{align}
  \label{Q0forP2}
  Q^0 \
  &= - \tfrac{1}{2} \big( G_s \!\cdot\! \Lop G + G_u \!\cdot\! \Pop G \big) ,
  \\[8pt]
  \label{Q1forP2}
  Q^1_{x_1 x_1} 
  &=	G_{v u} - 2G_{w s} - G_{wss}                     
	+ G_{s} G_{u s s}
	+ G_{u} \!\cdot\! \Lop G_{s s}
	+ \big( G_{v s} \!\cdot\! \Lop G_{s s}
	+ G_{v u} \!\cdot\! \Pop G_{s s} \big) ,
	\\[8pt]
	\label{Q2forP2}
		Q^2_{x_1 x_1} 
		&= 
	G_{wu} 	+ G_u G_{uss}
	+ \big( G_{ws} \!\cdot\! \Lop G_{ss}
	+ G_{wu} \!\cdot\! \Pop G_{ss} \big) .
\end{align}
(In deriving (\ref{Q1forP2}), the chain rule enters non-trivially, producing
five extra terms which are exactly minus the right hand side of
Equation~(\ref{eq2GP2}), which is then used backwards.)

\end{blanko}

\begin{blanko}{The locus of marked cusp.}
	In the space $M_{0,1}(\P^2,d)$ of irreducible maps with
	a single mark $p_1$,  Consider the locus of maps such that $p_1$ is a critical
	point, i.e.\ the differential vanishes at $p_1$.  The locus of non-immersions
	is of codimension 1, so  requiring further that the mark is critical gives
	codimension 2.  Let $\cusp_1$ denote the closure of this locus in $\ov 
	M_{0,1}(\P^2,d)$, the locus of maps having a cusp at $p_1$.  In
	spaces with more marks, $\cusp_1$ is defined as the pull-back of the one
	in $\ov M_{0,1}(\P^2,d)$ via the forgetful morphism.
\end{blanko}

\begin{prop}\label{cusp}
The class of this marked cusp locus is
$$
\cusp_1 = 3\evalclass_1^2 + 3\evalclass_1\mpsiclass_1 + \mpsiclass_1^2 
- \Pi_1 .
$$
\end{prop}

\begin{bevis}
  We start out with a family of stable un-pointed maps
$$\begin{array}{rlcc}
&\XX		& \overset{\displaystyle{\mu}}{\longrightarrow}		& \P^2		\\
\phantom{\pi}\pi\!\!\!\!\!&\downarrow&  		& 		\\
&B		& 		& 
\end{array}$$
where $B$ and $\XX$ are smooth, and the locus $N\subset \XX$ of singular points
of the fibres is of codimension 2.  Let $I\subset\OO_\XX$ be the ideal sheaf of
$N$. The exact sequence
$$
0 \pil \pi\upperstar \Omega_B \pil \Omega_{\XX} \pil I \tensor\kan_\pi\pil 0
$$
yields the relation of total Chern classes
$
c(T_{\XX}) = \pi\upperstar c(T_B) \big( 1 - K_\pi+ 
[N]\big),
$
and thus
\begin{equation}\label{-N}
\frac{\pi\upperstar c(T_B)}{c(T_\XX)}
= 1 + K_\pi + K_\pi^2 -  [N].
\end{equation}
Here, $K_\pi \df c_1(\kan_\pi)$, and we also set $H\df\mu\upperstar 
c_1(\OO(1))$.

Denote temporarily by $D$ the class of the locus of points in $\XX$ where the
differential $T_\XX \to (\pi\times\mu)\upperstar T_{B\times\P^2}$ fails to have
rank $2$.  By Porteous' formula, $D$ is the degree-$2$ part of the total
Chern class
$$
\mu\upperstar  c(T_{\P^2}) \cdot \frac{\pi\upperstar c(T_B)}{c(T_\XX)}
= (1+3H+3H^2)(1 + K_\pi + K_\pi^2 -  [N]),
$$
by~(\ref{-N}). In other words,
$$
D = 3H^2 + 3H K_\pi + K_\pi^2 - [N]
$$
All this is basically \S 4.d of Diaz-Harris~\cite{Diaz-Harris}.

Now equip the family with a section $\sigma_1:B\to\XX$ that transversely 
intersects $N$.  The marked-cusp class of the family is just
$\cusp_1 = \sigma_1\upperstar D$.   Now $\sigma_1\upperstar H = \evalclass_1$ and
$\sigma_1\upperstar K_\pi = \psiclass_1=\mpsiclass_1$, so we get
$$
\cusp_1 = 3\evalclass_1^2 + 3\evalclass_1\mpsiclass_1 + \mpsiclass_1^2 
- \sigma_1\upperstar [N].
$$
(In a family with more sections, we must pull back these constructions;
therefore the modified psi class is the correct one to use.)
This family of marked maps is not stable, but there is a well-defined
stabilisation.  It only remains to notice that the locus $\sigma_1\upperstar [N]
\subset B$ of the unstable family is the same as $\Pi_1$ of the stabilised
family.
\end{bevis}

\begin{blanko}{An alternative construction},\label{alternative}
  also given in \cite{Kock:thesis}, describes $\cusp_1$ as the locus of maps
  $\mu:C\to\P^2$ such that a whole pencil of lines in $\P^2$ are tangent to
  $\mu(C)$ at $\mu(p_1)$.  In other words, it is the degeneracy locus of the map
  of vector bundles $\sigma_1\upperstar V^3\to \sigma_1\upperstar L^2$, where
  $V^3$ is the $\mu$-pull-back of the complete linear system $H^0(\P^2,\OO(1))$,
  and $L^2$ is the sheaf of first principal parts of $\mu\upperstar \OO(1)$. 
  But then it is necessary to correct for $\Pi_1$ afterwards.
\end{blanko}

\begin{BM}
	For $d=1$, the locus is empty, so in this case Porteous' formula yields the
	relation $3\evalclass^2 + 3\evalclass_1\mpsiclass_1 + \mpsiclass_1^2 = 0$. 
	Under the natural identification of $\ov M_{0,1}(\P^2,1)$ with the incidence
	variety $I\subset \P^2\times\check \P^2$ of points and lines in $\P^2$, this
	relation is equivalent to the well-known relation $h^2 + \check h^2 =
	h\check h$.
	
	For $d=2$, the multiple-covers occur already in codimension 1.  On the
	other hand, there are no birational maps in degree 2 with a critical point,
	so for $d=2$ the locus $\cusp_1$ consists of all the double covers such that
	the mark is one of the ramification points.
	
	For $d\geq 3$, the locus of multiple-covers is of codimension at least 2, so
	the extra condition of having the mark as one of the ramification points
	prevents these curves from contributing.  So in this case the locus
	$\cusp_1$ consists generically of birational maps.
\end{BM}
\begin{blanko}{Further cusp conditions.}\label{cusp-at-Q}
	Consider the codimension-3 condition of the marked cusp mapping to a given
	line.  The class $\cuspl_1$ of this condition is obtained simply by cutting
	with $\evalclass_1$:
$$
\cuspl_1 = 3\evalclass_1^2\mpsiclass_1 + 
\evalclass_1\mpsiclass_1^2 - \evalclass_1\Pi_1.
$$
Similarly, the locus of maps with marked cusp mapping to a specified point is
$$
\cuspp_1 = \evalclass_1^2\mpsiclass_1^2 - \evalclass_1^2\Pi_1.
$$
These two loci can also be constructed by the approach of \ref{alternative},
replacing the complete linear system by smaller systems, cf.~\cite{Kock:thesis}.
\end{blanko}

\begin{blanko}{Tangency conditions in cuspidal environment.}
	Suppose we are inside the locus $\cusp_1$ and want to impose the condition of
	being tangent to a given line $L$ at another mark, say $p_2$.  Since for the
	general map in $\cusp_1$, the differential vanishes only simply a $p_1$, the
	arguments of [GKP] 3.1 and 3.3 show that the locus of maps which are not
	transversal to $L$ at $p_2$ is reduced of class $\evalclass_2 (\evalclass_2 +
	\mpsiclass_2)$.  However, contrary to the case of nodal curves, this locus
	has two irreducible components.  In addition to the locus of honest
	tangencies, there is a component consisting of maps such that the $p_1$-cusp
	maps to $L$ and the two marks have come together, i.e., $
	\cusp_1\!\cdot\!\evalclass_1 \!\cdot\!\deltaclass_{12}$.  We do not want to
	count these maps as tangencies, so in conclusion, the class of $p_2$-tangency
	in $p_1$-cuspidal environment is
\begin{equation}\label{Theta'}
\Theta_2 ' = \evalclass_2 (\evalclass_2 + \mpsiclass_2) - 
\evalclass_1\deltaclass_{12}.
\end{equation}

Similarly, in $p_1$-cuspidal environment the class of $p_2$-tangency to a given
line at a specified point is
$$
\Xi_2' = \evalclass_2^2 \mpsiclass_2 - \evalclass_1^2\deltaclass_{12}.
$$

We can now apply these conditions iteratively, and the top intersections will
be the characteristic numbers
for cuspidal plane curves.
\end{blanko}

Using the generating functions for the slightly enriched first enumerative
descendants, and their relation to the tangency potential, it is straightforward
to derive differential equations determining the cusp characteristic numbers
from the nodal ones.
Let $C_d(a,b,c)$ denote the number of cuspidal plane curves passing through $a$
points, tangent to $b$ lines, and tangent to $c$ lines at specified points.  Let
$\CL_d(a,b,c)$ be defined similarly but requiring the cusp to fall on a specified
line, and let $\CP_d(a,b,c)$ denote the numbers where the cusp is required to
fall at a specified point.  Let $K(s,u,v,w)$, $\KL(s,u,v,w)$ and $\KP(s,u,v,w)$ be
the corresponding generating functions (the formal variables being defined as in
\ref{setupG}).

\begin{prop}\label{cuspDE}
  The cusp potentials $\KP$, $\KL$, and $K$ are determined by the characteristic
  number potential $G$ through the following equations.
\begin{align}
  \label{KP}
  	\KP_{s s} \ &= \ G_{w u} - G_{u s}G_{u s}
	+ \big( G_{w s} \!\cdot\! \Lop G_{ss}
	+ G_{w u} \!\cdot\! \Pop G_{ss} \big) ,
	\\[8pt]
	\label{KL}
		\KL_{s s} \ &= \ + G_{v u} + 2 G_{w s s}  - v \KP_{s s}
		- 2 G_{w s}  - G_{us}G_{ss} - G_{us} \!\cdot\! \Lop  G_s
		\\ \notag & \phantom{xxxxxxxxxxxxxxxxxxxxxx}
	+ \big( G_{v s} \!\cdot\! \Lop G_{ss}
	+ G_{v u} \!\cdot\! \Pop G_{ss}\big)  
	 ,
	\\[8pt]
	\label{K}
		K \ &= \ 3 G_{v} 
	- v \KL - (\brok{1}{2}v^2 + w) \KP 
	- \big( G_s \!\cdot\! \Lop G + G_u \!\cdot\! \Pop G \big)
	.
\end{align}
\end{prop}
\begin{bevis}
  The main point is to eliminate the diagonal classes.
  In each term of the
  expansion of the top product, the diagonal class $\deltaclass_{1i}$ is alone
  at mark $p_i$, so we can push down forgetting $p_i$.  The push-down formula is
  simply $\pi_i{}\lowerstar \deltaclass_{1i} = 1$ (cf.~[GKP], 1.3.2.).
  
  Since $\evalclass_1^3=0$, and since all diagonal classes come accompanied by a
  factor $\evalclass_1$, only few diagonal class terms survive the expansion of
  the top product.  In the presence of a factor $\cuspp_1 =
  \evalclass_1^2(\mpsiclass_1^2-\Pi_1)$, all the diagonal classes of the top
  product vanish.  Thus,
  $$
  \KP ( s, u, v, w) = (Q^2-P^2)(x_1, x_2, y_1, y_2).
  $$
  Now take double derivative with respect to $s=x_1$ and apply
  Equations~(\ref{Q2forP2}) and (\ref{P2forP2}).  This establishes (\ref{KP}).
	
  In the integral corresponding to (\ref{KL}), since there is a factor
  $\evalclass_1$ in $\cuspl_1$, there is room for at most one diagonal class
  in each term of the expansion.  So we get
\begin{eqnarray*}
	d^2 \; \CL_d(a,b,c) & = & d^2\; \cuspl \; \Omega^a \; \Theta'{}^b \; \Xi^c \\
	 & = & d^2\; \cuspl \; \Omega^a \; \Theta^b \; \Xi^c  - d^2\, b \,\cuspp \; 
	 \Omega^a \; \Theta^{b-1} \; \Xi^c .
\end{eqnarray*}
Here $\cuspp_1$ arises as $\evalclass_1\cdot \cuspl_1$.
The last term explains $- v \KP_{s s}$ in the formula.  In the first term 
we plug in $\cuspl_1= 3\evalclass_1^2\mpsiclass_1 + 
\evalclass_1\mpsiclass_1^2 - \evalclass_1\Pi_1 = 3\Xi_1 + 
\evalclass_1(\mpsiclass_1^2 - \Pi_1)$. Thus
$$
\KL_{s s} = - v \KP_{s s} + 3 G_{wss} + (Q^1_{x_1 x_1} - P^1_{x_1 x_1}) .
$$
The result now follows from Equations~(\ref{P1forP2}) and (\ref{Q1forP2}).

  Finally in the expansion of the integral corresponding to (\ref{K}), we get
  $b$ terms corresponding to one diagonal class from $\Theta'$, further
  $\binom{b}{2}$ terms with two diagonal classes from $\Theta'$, and finally $c$
  terms with one diagonal class from $\Xi'$:
\begin{eqnarray*}
	C_d(a,b,c) & = & \cusp  \; \Omega^a \; \Theta'{}^b \; \Xi'{}^c \\
	 & = & \cusp  \; \Omega^a \; \Theta^b \; \Xi^c
	 -b \cuspl  \; \Omega^a \; \Theta^{b-1} \; \Xi^c
	 +\smallbinom{b}{2} \cuspp  \; \Omega^a \; \Theta^{b-2} \; \Xi^c
	 -c \cuspp  \; \Omega^a \; \Theta^b \; \Xi^{c-1}   \\
&=& \cusp  \; \Omega^a \; \Theta^b \; \Xi^c 
	   -b \CL_d(a,b-1,c) -\smallbinom{b}{2} \CP_d(a, b-2,c) - c \CP_d(a,b,c-1).
\end{eqnarray*}
The last three terms explain $- v \KL - (\brok{1}{2}v^2 + w) \KP$ in  the formula.
The first term is expanded to $$
\cusp  \; \Omega^a \; \Theta^b \; \Xi^c =
3 N_d(a,b+1,c) + (\mpsiclass_1^2 -\Pi_1) \; \Omega^a \; \Theta^b \; \Xi^c,
$$
and this last term corresponds to $Q^0- P^0$ which is then expanded using 
Lemma~\ref{P0+Q0} and Equation~(\ref{Q0forP2}).
\end{bevis}

These differential equations are very similar to the recursions used in
Ernstr\"om-Kennedy~\cite{EK1} (found with completely different methods), and are
presumably equivalent (modulo Equations (\ref{eq1GP2}) and (\ref{eq2GP2})), but
I have not been able to identify all the terms of their recursion.

\begin{BM}
	Setting $v = w = 0$ (corresponding to considering only incidence
	conditions) and then differentiating with respect to $s$ yields
	\begin{eqnarray*}
		K_{s} & = & 3 G_{vs} - \brok{\d}{\d s}G_{s}^2  \pmod{ (v,w)}\\
		 & = & 3 \big( G_{us} - G_{u} +\brok{1}{2}G_{s s}^2\big)
	- \brok{\d}{\d s}G_{s}^2 \pmod{ (v,w)} ,
	\end{eqnarray*}
	which is equivalent to the recursion of Proposition~5 in 
	Pandharipande~\cite{Pand:9504}.
\end{BM}

\section{Cuspidal curves in $\P^1\times\P^1$}
\label{Sec:P1xP1}

\begin{blanko}{Set-up for $\P^1\times\P^1$.}\label{P1xP1}
	Let $T_0$ be the fundamental class; let $T_3$ be the class of a point; and
	let $T_1$ and $T_2$ be the hyperplane classes pulled back from the two
	factors.  A curve of class $\beta$ is said to have bi-degree $(d_1,d_2)$,
	where $d_1=\int_\beta T_1$ and $d_2=\int_\beta T_2$.  A curve of bi-degree 
	$(1,0)$ is called a horizontal rule, and a curve of bi-degree $(0,1)$ a 
	vertical rule. 
	
	Let $N_{(d_1,d_2)}(a,b,c)$ denote the characteristic numbers of irreducible
	rational curves in $\P^1\times\P^1$ of bi-degree $(d_1,d_2)$ passing through
	$a$ general points, tangent to $b$ general curves of
	bi-degree $(1,1)$, and tangent to $c$ such curves at
	specified point.  
	The classes corresponding to these three conditions are, respectively:
	$\Omega=\mtau_0(T_3)$, $\Theta=2\mtau_0(T_3) + \mtau_1(T_1) + \mtau_1(T_2)$,
	and $\Xi=\mtau_1(T_3)$.

Let $G(u_1,u_2,u,v,w)$ be the
	corresponding generating function ($u_1$ and $u_2$ being the formal variables
	corresponding to the partial degrees $d_1$ and $d_2$).
	Then we have $G(u_1,u_2,u,v,w) =
\Gamma(x_1,x_2,x_3,y_1,y_2,y_3)$, with $x_1 = u_1$, $x_2=u_2$, $x_3= u +
2 v$; $y_1=v$, $y_2=v$, $y_3=w$.  For convenience, put also $$s\df u_1 + u_2,$$
the formal variable corresponding to $T_1 + T_2$.  We have
$$
(\gamma^{ef}) = 
\begin{pmatrix}
0 &      0 &      0 &            1        \\
0 &     0  &          1 &     2 y_1        \\
0 &      1 &    0 &    2 y_2        \\
1 &2 y_1&2 y_2    & 4 y_1 y_2 + 2 y_3
\end{pmatrix}
=
\begin{pmatrix}
0 &      0 &      0 &            1   \\
0 &     0  &          1 &     2 v    \\
0 &      1 &    0 &    2 v        \\
1 &2 v&2 v    & 4 v^2 + 2w      
\end{pmatrix} .
$$
Define three differential operators corresponding to the three last lines of 
this matrix,
\begin{align*}
\Lopet &\;\df\; \phantom{2v}\brok{\d}{\d u_2} + 2v \brok{\d}{\d u} ,\\
\Lopto &\;\df\; \phantom{2v}\brok{\d}{\d u_1} + 2v \brok{\d}{\d u} ,\\
\Pop &\;\df\; 2v\brok{\d}{\d u_1} + 2v\brok{\d}{\d u_2} + (4v^2+2w) 
\brok{\d}{\d u} ,\\
\intertext{and for convenience put also}
\Lop &\;\df\; \Lopet+\Lopto = \brok{\d}{\d s} + 4v \brok{\d}{\d u}
.
\end{align*}

Equations (25) and (26) of [GKP] read
\begin{align}
  \label{eq1GP1xP1}
G_{vs} &= 2 G_{us} - 2 G_u + \tfrac{1}{2}\big( G_{s u_1} \!\cdot\!  \Lopet G_{s}
+G_{s u_2} \!\cdot\!  \Lopto G_{s} + G_{us} \!\cdot\! \Pop G_{s} \big)  ,
\\[8pt]
\label{eq2GP1xP1}
G_{wss} &=
2 G_{uu} + \big( G_{u u_1} \!\cdot\! \Lopet G_{ss} 
+G_{u u_2} \!\cdot\! \Lopto G_{ss}
+ G_{uu} \!\cdot\! \Pop G_{ss} \big).
\end{align}
%
%
\end{blanko}


\begin{blanko}{Differential equations for the slightly enriched potentials.}
  \label{PQforP1xP1}
  We have
	$$
(\gamma_{(12)}^{ef}) = 
\begin{pmatrix}
0 &      0 &      0 &            0       \\
0 &     0  &          0 &    1        \\
0 &      0 &    0 &    1        \\
0 &1&1   & 2 y_1 + 2 y_2      
\end{pmatrix}
\qquad\text{ and }\qquad
(\gamma_{3}^{ef}) = 
\begin{pmatrix}
0 &      0 &      0 &            0       \\
0 &     0  &          0 &    0        \\
0 &      0 &    0 &    0        \\
0 &0  &0   & 1      
\end{pmatrix} .
$$
Now applying the coordinate changes to Propositions~\ref{PkGamma} and
\ref{QkGamma} we get:
\begin{align}
\label{P12forP1xP1}
P^{(12)}_{x_1 x_1} 
&= G_{us} G_{ss} + G_{uss} G_s + G_u \!\cdot\! \Lop G_{ss} 
+ G_{us} \!\cdot\! \Lop G_s,
\\[8pt]
\label{P3forP1xP1}
P^3_{x_1 x_1} 
&= G_{us}G_{us} + G_{uss} G_u ,
\\[8pt]
\label{Q0forP1xP1}
Q^0 \phantom{xx} &=
-\tfrac{1}{2} \big( G_{u_1}\!\cdot\!\Lopet G + G_{u_2}\!\cdot\!\Lopto G
+ G_{u}\!\cdot\!\Pop G \big) ,
\\[8pt]
\label{Q12forP1xP1}
Q^{(12)}_{x_{12}x_{12}} 
&= 
	 2 G_{v u} - 4 G_{w s} - 2 G_{ w s s} + G_{s}G_{u s s}
		+ G_{u}\!\cdot\!\Lop G_{ss} \\
		\notag &\phantom{xxxxxxxxxxxxxx}
+ \big( G_{v u_1}\!\cdot\!\Lopet G_{ss} + G_{vu_2}\!\cdot\!\Lopto G_{ss}
+ G_{vu}\!\cdot\!\Pop G_{ss} \big) .
\\[8pt]
\label{Q3forP1xP1}
Q^3_{x_{12}x_{12}} 
&= 		2 G_{wu} + G_{u}G_{u s s}
+	\big( G_{w u_1}\!\cdot\!\Lopet G_{ss} + G_{wu_2}\!\cdot\!\Lopto G_{ss}
+ G_{wu}\!\cdot\!\Pop G_{ss} \big)
		 .  \phantom{xx}
\end{align}
The derivation of these formulae follows the same arguments as in \ref{PQforP2}.

\end{blanko}
%

\begin{blanko}{Differential equations.}
  Deriving equations for the cusp potentials for $\P^1\times\P^1$ similar to
  those of \ref{cuspDE} is now straightforward.  Since the
  tangent bundle has total Chern class $1 + 2(T_1+T_2) + 4T_3$, the locus of
  cusp at mark $p_1$ is
$$
4 \nu_1\upperstar (T_3) + 2(\nu_1\upperstar (T_1)+\nu_1\upperstar 
(T_2))\mpsiclass_1 + \mpsiclass_1^2 - \Pi_1.
$$
Let $\KP$ be the potential corresponding to cusp mapping to a specified point,
(and further $a$ conditions of passing through a point, $b$ conditions of being
tangent to a $(1,1)$-curve, and $c$ conditions of tangenciating such a curve at
a specified point).  Then
\begin{align*}
\KP_{ss} &= 2 G_{wu} - G_{us}G_{us}
+ \big( G_{w u_1}\!\cdot\!\Lopet G_{ss} + G_{wu_2}\!\cdot\!\Lopto G_{ss}
+ G_{wu}\!\cdot\!\Pop G_{ss} \big)
		 .   
\end{align*}
Let $\KL$ be the generating function for such characteristic numbers, but with 
the cusp mapping to a specified $(1,1)$-curve. Then
\begin{align*}
\KL_{s s} 
&= 2 G_{v u} -4 G_{w s} - 2v \KP_{s s}  
+ 2 G_{w s s} - G_{u s}G_{s s} - G_{us}\!\cdot\!\Lop G_s \\
& \phantom{xxxxxxxxxxxx}
+ \big( G_{v u_1}\!\cdot\!\Lopet G_{ss} + G_{vu_2}\!\cdot\!\Lopto G_{ss}
+ G_{vu}\!\cdot\!\Pop G_{ss} \big) .
\end{align*}
And finally, let $K$ be the generating function for the characteristic numbers
of cuspidal curves in $\P^1\times\P^1$, with the cusp varying freely. Then
\label{KP1}
\begin{align*}
	K &= 2G_{v} - v \KL - (v^2 + w) \KP 
	-  \big( G_{u_1}\!\cdot\!\Lopet G + G_{u_2}\!\cdot\!\Lopto G
+ G_{u}\!\cdot\!\Pop G \big)
.
\end{align*}
\end{blanko}
\begin{blanko}{Enumerative significance.}\label{cuspPP}
  A priori these numbers count
  also reducible curves, one of whose twigs is a multiple cover of a rule.  In
  fact, already the locus $\cusp_1$ is not irreducible: it has a component for
  each boundary divisor corresponding to degree splitting $(m,n) = (i,0) +
  (m-i,n)$.  For each of these divisors, the one-primed twig is always a
  multiple-cover of a horizontal rule, and forcing the mark to a
  ramification point produces the `cusp' already in codimension 2.  The other
  ramification points can then satisfy tangency conditions, giving contribution
  in the characteristic number.  (Similarly of course for maps comprising a
  cover of a vertical rule.)
  
  However, when there are no conditions on the cusp, all solutions are in fact
  irreducible curves. 
  This happens because one degree of freedom (that of varying the position of
  the ramification point marked $p_1$ which counts as the cusp), is useless for
  the sake of satisfying tangency (or incidence) conditions, since we have
  already excluded the case where the tangency condition is fulfilled at $p_1$. 
  Now the dimension count is easy: The multiple-cover twig has $2i-2$
  ramification points, of which $2i-3$ can be used for satisfying tangency
  conditions.  The twig can satisfy a single incidence condition (which
  completely fixes the support of the curve).  The honest twig (of bi-degree
  $(m-i,n)$) can as usual satisfy $2m-2i+2n-1$ conditions.  Thus we get a total
  of $2m+2n-3$ degrees of freedom, while the number of conditions imposed is
  $2m+2n-2$.  So no reducible curves can contribute.
  
\bigskip

  The situation is different in the cases where the cusp is required to
  fall on a $(1,1)$-curve or at a point.  In these cases, the cusp can account
  for one condition in addition to just being a cusp.  For example, in the case
  where the cusp is required to fall on a given point: The given point fixes
  the rule supporting the image of the $p_1$-twig, and this must be a
  ramification point.  Then further $2i-3$ ramification points can account for
  as many tangency conditions.  On the other twig there is room for $2m-2i+2n
  -1$ conditions.  Total: $2m+2n-4$, which is exactly the number of conditions
  imposed.  So the potential $\KP$ encodes also reducible curves.
  A similar observation applies to $\KL$.
  
  \medskip
  
  Correcting for these reducible curves is a case for the techniques described
  in [GKP] \S 5.  Since the unwanted
  contribution are multiple covers, the correction terms involve the (genus
  zero) Hurwitz numbers.
\end{blanko}

\begin{blanko}{Hurwitz numbers and multiple-covers of a rule.} 
	(Cf.~[GKP] \S 5.)  Consider
	$X=\P^1$ (with $T_0$ = fundamental class, $T_1$ = class of a point).  The
	invariants $N_d(b) = \brac{ \mtau_1(T_1)^{b} }\!_{d}^{\P^1}$ are the simple
	genus zero Hurwitz numbers (the number of $d$-sheeted coverings $\P^1\to
	\P^1$ simply ramified over $b=2d-2$ given points).  The corresponding
	potential,
	\begin{equation}\label{Hurwitz}
		H(t,v) =
	\sum_{d>0} \exp(dt) \sum_b \frac{v^b}{b!} N_d(b) .
	\end{equation}
	satisfies the topological recursion relation
	$H_{vt} = v H_{tt} \!\cdot\! H_{tt}$.

	For $\P^1\times\P^1$, the generating function for the covers of a horizontal
	rule is
	\begin{equation}
	I(u_1,u,v,w) = uH_{u_1} + (v^2+w)H_v ,\label{rulecovers}
	\end{equation}
	where $H=H(u_1,v)$ is the Hurwitz potential of (\ref{Hurwitz}).  Indeed, the
	supporting rule for an $i$-sheeted map is fixed by either one incidence
	condition, two tangency conditions, or one flag condition.  Once the
	supporting rule is fixed, the Hurwitz potential encodes the number of
	possible coverings.  For the incidence condition, there are $i$ choices for
	the mark; this explains the factor $uH_{u_1}$.  For the case of $2i-1$
	tangency conditions, the rule must pass through one of the two intersection
	point of two of the given curves.  This gives $2\!\cdot\!\binom{2i-1}{2}$
	choices for the supporting rule, explaining the term
	$2\!\cdot\!\frac{v^2}{2!} H_v$.  Finally, for one flag condition and $2i-3$
	tangency conditions, the flag fixes the supporting rule and translates into
	an extra $v$ condition on the covering of that rule.
		
	For coverings of a vertical rule, we similarly find the generating function
	\begin{equation}
		J(u_2,u,v,w) = uH_{u_2} + (v^2+w)H_v ,
	\end{equation}
	with $H=H(u_2,v)$.
\end{blanko}
Now the correction term corresponding to the fake reducible cusps counted in $\KP$ 
is described by:
\newcommand{\irr}{\operatorname{irr}}
\begin{prop}
  The generating functions $K^{\irr}$, $\KL^{\irr}$, and $\KP^{\irr}$ for the
  characteristic numbers of irreducible cuspidal curves in $\P^1\times\P^1$ are
  related to $K$, $\KL$, and $\KP$ like this:
  \begin{align*}
	 K^{\irr} \ \ &= \ K  ,
	 \\
    \KL^{\irr} \ &= \ \KL - 
   \big( I_{v u_1} \!\cdot\! \Lopet G 
	+ J_{v u_2} \!\cdot\! \Lopto G 
	+ (I_{vu} \! + \! J_{vu}) \!\cdot\! \Pop G \big) ,
	\\
	\KP^{\irr} \ &= \ \KP - 
   \big( I_{w u_1} \!\cdot\! \Lopet G  
  + J_{w u_2} \!\cdot\! \Lopto G   \big) .
  \end{align*}
\end{prop}
\begin{bevis}
  It has already been shown that $K^{\irr} =  K$.  For the others,  
  the correction term has an $I$-part coming from multiple-covers of the
  horizontal rule and a $J$-part corresponding to covers of a vertical rule. 
  Let us explain the $I$-part.  The arguments of [GKP] \S 5 show
  that
  \begin{equation}
	 I_{u_1} \!\cdot\! \Lop{}_1 G + I_u \!\cdot\! \Pop G
	 \label{IG}
  \end{equation}
  is the generating function for the numbers of the reducible maps comprising a
  multiple-cover of a horizontal rule.  In the present situation $p_1$ is a
  ramification point on the twig covering the horizontal rule.  Therefore,
  requiring further that $p_1$ maps to a given line is equivalent to requiring
  tangency to that line at $p_1$, so the $I$-potentials in (\ref{IG}) acquire a
  derivative with respect to $v$.  This accounts for the $I$-part of the
  correction term to $\KL$.  Concerning $\KP$: requiring the ramification point
  $p_1$ to map to a given point is equivalent to imposing a flag condition on
  $p_1$, so in this case the $I$-potentials in (\ref{IG}) acquire a derivative
  with respect to $w$.  It remains to note that $I_{wu}=0$ since a
  multiple-cover of a rule cannot pass through two given general points.
\end{bevis}

\begin{eks}\label{Kinc} 
  Setting $v=w=0$ in the equation for $K$ yields an easy expression for the
  numbers $C_{(m,n)}$ of (irreducible) cuspidal curves in $\P^1\times\P^1$ of
  bi-degree $(m,n)$ that pass through $2m+2n-2$ general points, in terms of the 
  corresponding numbers $N_{(m,n)}$ of nodal curves:
$$
C_{(m,n)} = \frac{4(d-1)}{d} N_{(m,n)} + \sum \binom{2d-2}{2d'-1} 
\frac{(m' n'' + n' m'')(d' d'' - d )}{d} N_{(m',n')}N_{(m'',n'')} ,
$$
where the sum is over degree splittings $m'+m''=m$ and $n'+n''=n$.
For short we have set $d=m+n$, $d'=m'+n'$, and $d''=m''+n''$.
\end{eks}

\section[Characteristic numbers of curves with triple contact]%
{Characteristic numbers of curves with \\ a prescribed triple contact}

\label{Sec:flex}

Let $V\subset X$ be a
general, very ample hypersurface, given as the zero scheme of a section $f$ of a
line bundle $L$, and set $\evalclass_i \df c_1(\nu_i\upperstar L)$.

\begin{blanko}{Components mapping into $V$.}
	Denote by $\lineartwiginH_1= \lineartwiginH_1(V)$ (the closure of) the locus
	of maps such that $p_1$ is on a non-contracting twig that maps entirely into
	$V$.  Note that we define $\lineartwiginH_1$ only on set-theoretic level and
	not as a cycle.  In general, this locus has irreducible components in various
	codimensions.  The locus of {\em irreducible} such maps
	$\lineartwiginH_1\subset M_{0,n}(X,\beta)$ is the zero scheme of a regular
	section of the vector bundle $\pi\lowerstar \mu\upperstar L$, so the
	codimension is $\dim H^0(\P^1, \mu\upperstar L)$ (or it is empty).
	
	However, commonly the component consisting of irreducible such maps is not
	the one of lowest codimension: For example, if $X=\P^r$ and $V$ is a
	hypersurface of degree $z$, then the locus $\lineartwiginH_i\subset
	M_{0,n}(\P^r,d)$ is of codimension $dz+1$, while the locus of reducible such
	maps having $p_i$ on a twig of degree $1$ is of codimension $z+2$.  So for
	$d\geq 2$ or $z\geq 2$, the locus $\lineartwiginH_1$ is always of codimension
	at least 3.  The only situation in which we need to worry about codimension
	2, is $d=z=1$: in this case $\codim\lineartwiginH_1 =2$, and of course these
	maps are irreducible.
\end{blanko}
\begin{blanko}{Enumerative irrelevance of $\lineartwiginH_1$.}
	Note first of all that $\lineartwiginH_1$ is compatible with inverse image
	under forgetful morphisms, since in the definition we have excluded the case
	of $p_1$ belonging to a twig contracting to a point in $H$.
	
	Now in the locus $\lineartwiginH_1$, the mark $p_1$ has not been fixed: it
	can be any point of the twig mapping into $V$, and since this twig is not
	just a contracting twig, different choices of where to put the mark on the
	twig are non-isomorphic.  The consequence of this observation is that
	$\lineartwiginH_1$ has zero push-down under forgetting $p_1$ (whatever the
	multiplicities attributed to the components of it).  This is turn implies
	the vanishing of any top product where
	$\lineartwiginH_1$ is up against pull-back-compatible classes belonging to
	other marks.
\end{blanko}

With these remarks we are in position to describe the locus of maps having a 
triple contact with $V$.
For the argument to work we must impose the following condition on the
moduli space and on $V$:

\begin{list}{}{%
	\setlength{\topsep}{9pt}
	\setlength{\leftmargin}{44pt}
	\setlength{\rightmargin}{10pt}%
	\setlength{\labelwidth}{20pt} \setlength{\itemsep}{3pt}
	\setlength{\parsep}{0pt} } 
	
	\item[{(\(\star\))}] The general map tangent to $V$ is only simply tangent to $V$.
\end{list}
The easiest example violating this condition is $\ov M_{0,1}(\P^r,1)$ when
$V$ is a hyperplane: every line tangent to $V$ is
actually contained in $V$.

\begin{prop}\label{flex}
	With (\(\star\)) satisfied, let $\flex_1$ denote the closure of the locus of
	irreducible immersions having a triple contact to $V$ at the mark $p_1$. 
	Then its class is
$$
\flex_1 =
\evalclass_1(\evalclass_1+\mpsiclass_1)(\evalclass_1+2\mpsiclass_1)
	- [\lineartwiginH_1]_3 - \evalclass_1\Pi_1,
	$$
	where $[\lineartwiginH_1]_3$ denotes the codimension-3 part of
	$\lineartwiginH_1$.
\end{prop}

\begin{bevis}
	We can assume there is only one mark $p_1$.  Consider the zero
	scheme of the section $\sigma_1\upperstar \d_\pi^2\mu\upperstar f$ of the
	bundle of second jets, $\sigma_1\upperstar J_\pi^2\mu\upperstar L$.  It 
	follows from the standard jet short exact sequence that the class of this 
	locus is
	$$
	c_3(\sigma_1\upperstar J_\pi^2\mu\upperstar L) 
	= \evalclass_1(\evalclass_1+\mpsiclass_1)(\evalclass_1+2\mpsiclass_1).
	$$
	It is the
	codimension-3 locus of maps such that the pull-back of $f$ vanishes to order
	3 at $p_1$.  If it happens isolated at $p_1$, then we are in $\flex_1$. 
	Note that there is no contribution from $p_1$-marked cusps: the locus of such 
	cusps mapping to $V$ is of codimension 3, but in general the
	differential vanishes only simply and thus the first jet vanishes but not the 
	second.
	
	If $\mu\upperstar f$ vanishes identically on the twig carrying $p_1$ then
	{\em either} this twig is of degree zero, and then we are at $\evalclass_1
	\Pi_1$ (or in higher codimension); {\em or} the twig is of positive degree,
	and that means it maps entirely inside $V$ and then we are in the locus
	$\lineartwiginH_1$.  (By the condition (\(\star\)), no component of
	$\lineartwiginH_1$ can have codimension lower than 3.  There may be
	components of higher codimension but since a section of a rank-3 vector
	bundle cannot vanish isolated along such a component, we conclude that these
	loci are in fact already contained in one of the codimension-3 loci, which
	must then be $\flex_1$.  In other words, such maps are limits of honest
	triple contacts.)
	
	To establish multiplicity-1 along $\Pi_1$ it is enough to construct a
	$1$-parameter family running inside the locus of maps tangent to $V$ at
	$p_1$, and check that the section of the second jet bundle vanishes with
	multiplicity $1$ along the $\Pi_1$ of this family.  For simplicity we perform
	the argument only in the case $X=\P^2$.

Let $B\subset \A^1$ be a small open neighbourhood of $0$, and consider the
rational map
\begin{eqnarray*}
	\mu \ : \ B\times \P^1 & \langpil & \P^2  \\
	\big(b,[s:t]\big) & \longmapsto & \begin{bmatrix} 
	b^d \; f  \\
b (t-s)^2 \big( 	b^{d-1} \; g \; + s^{m-1} t^{n-1} \big) \\
	s^m t^n + b^m t^d + b^n s^d\end{bmatrix}.
\end{eqnarray*}%
where $m+n=d$.  As to the other symbols: $f$ is a homogeneous polynomial in
$s,t$ of degree $d$, and $g$ is homogeneous of degree $d-2$, with roots distinct
from $(t-s)$.

Equip the map with the constant section
	$$
	\sigma_1(b) \df (b,[1:1]).
	$$
The map has two base points in the central fibre $C$ (equation $(b=0)$), namely
$[0:1]$ and $[1:0]$.  (The rest of $C$ is contracted to the point $[0:0:1]$ in
$\P^2$.)  The map has been constructed such that a single blow-up resolves the
map.  The result is a family over $B$ of $1$-pointed stable maps of degree $d$,
whose central fibre is an element in $\Pi_1$ (with degree splitting $m+n=d$). 
Since we are interested only in what happens at the section, we can forget about
the blow-up.

Let $V\subset \P^2_{[x:y:z]}$ be the line with equation $y=0$,  then clearly all
members of the family are tangent to $V$ at the first mark.  More precisely,
$\mu\upperstar V= 2\sigma_1 + C + R$, where $R$ is some
residual curve not passing through the point $(0,[1:1])$.  The coefficient $1$ 
in front of $C$ is the exponent of $b$ in the second coordinate function.
The vanishing of the
second jet is read off in the family as the intersection of $\sigma_1$ with
$\mu\upperstar V- 2\sigma_1 = C + R$.  This is just the transverse intersection
at $(1,[1:1])$.  Thus, $\Pi_1$ has multiplicity $1$.
\end{bevis}

\begin{BM}
	Note that we are not excluding the possibility that all maps in $\flex_1$ are
	in fact contained in $V$.  This happens for example in $\ov M_{0,1}(\P^3,1)$
	when $V\subset \P^3$ is a quadric surface.
\end{BM}

Since the general map in $\flex_1$ is not a cusp, we can now impose further
tangency conditions without necessity of extra corrections. Thus we are in
position to compute the characteristic numbers of rational curves having a
triple contact with $V$.

\begin{eks}\textsc{The projective plane.}\label{flexinP2}
	Fix a general curve $V\subset \P^2$ of degree $z$.  Let $\NF^z_d(a,b,c)$
	denote the number of rational curves of degree $d\geq 2$ that have triple
	contact with $V$, are tangent to $b$ lines, tangent to $c$ lines at specified
	points, and furthermore pass through $a$ specified points, ($a+b+2c = 2d-3$). 
	Then
$$
\NF^z_d(a,b,c) = \int_{\ov M} \flex_1 \; \Omega^a \; \Theta^b \; \Xi^c .
$$
When expanding the class of $\flex_1$, we can ignore $\lineartwiginH_1$ since
it's push-down is zero, so we find
$$=
3z^2 \Xi \cdot \; \Omega^a \; \Theta^b \; \Xi^c
+ 2z\evalclass_1 \mpsiclass_1^2 \; \Omega^a \; \Theta^b \; \Xi^c
- z\evalclass_1 \Pi_1 \; \Omega^a \; \Theta^b \; \Xi^c.
$$
If we let $F^z$ be the generating function for these numbers we see that
$$
F^z( s, u, v, w) = 3z^2 G_{w}( s, u, v, w) + 2z Q^1(
\fat{x},\fat{y}) - z P^1( \fat{x},\fat{y})
$$
subject to the usual change of variables, cf.~\ref{setupG}.  Now take double
derivative with respect to $s$ and use the expressions for $P^1$ and $Q^1$
(given in Equations~(\ref{Q1forP2}) and (\ref{P1forP2})):
\begin{align*}
F^z_{s s} \ &= \ z\bigg( (3z-2) G_{wss} +2 G_{v u}
	-4 G_{w s}  + G_s G_{uss}
	+ G_u \!\cdot\! \Lop G_{ss} - G_{us}G_{ss} - G_{us} \!\cdot\! \Lop G_s
	\\[-13pt] &\phantom{xxxxxxxxxxxxxxxxxxxxxxxxxxxxxxxxxxxxxxxx}
	+2\big( G_{vs} \!\cdot\! \Lop G_{ss} + G_{vu}\!\cdot\! \Pop G_{ss}\big)
	 \bigg).
\end{align*}

Perhaps the most interesting case is when $z=1$ so we are talking about curves
with a specified flex line.  
\end{eks}

\begin{BM}
	Exploring the finer geometry of $\P^2$, through the space of stable lifts to
	its second Semple bundle, Colley-Ernstr\"om-Kennedy~\cite{CEK} have recently
	given a much more comprehensive formula.  Not only does their formula allow
	any number of triple contact conditions to given curves, but the given curves
	are also allowed to have double points and cusps!  The only drawback of their
	formula is that the given curves are not allowed to be lines or to have
	linear components.
\end{BM}

\begin{eks}\textsc{The quadric surface}\label{flexinPP}
	$X=\P^1\times\P^1$.
Continuing the notation from Section~\ref{Sec:P1xP1}, 
let $\NF_{(m,n)}(a,b,c)$ be the number of rational curves of bi-degree $(m,n)$
which have triple contact with a given curve of bi-degree $(1,1)$, are
tangent to $b$ other curves of bi-degree $(1,1)$, tangent to $c$ such curves at
specified points, and furthermore passing through $a$ specified points, ($a+b+2c
= 2m+2n-3$).

Let $F$ be the generating function for these characteristic numbers.  Then the
arguments  of the preceding example show that
$$
F = 6 G_{w}
+ 2 Q^{(12)}  - P^{(12)}
$$
subject to the usual variable changes $x_1 = u_1$, $x_2= u_2$, $x_3 = u +2v$,
$y_1=v$, $y_2=v$, and $y_3 = w$.  Taking double derivative with respect to
$s=u_1 + u_2 = x_1 + x_2$, and using Equations~(\ref{Q12forP1xP1}) and
(\ref{P12forP1xP1}) we get the formula
\begin{align*}
  F_{ss} \ &= \ 
  2 G_{wss} + 4 G_{v u} - 8G_{w s} + G_s G_{uss} + G_u \!\cdot\! \Lop G_{ss}
  - G_{us}G_{ss} - G_{us}\!\cdot\!\Lop G_s
  \\
& \phantom{xxxxxxxxxxxxxxxxxxxxxxxxxxx}
+ 2\big( G_{v u_1}\!\cdot\!\Lopet G_{ss} + G_{vu_2}\!\cdot\!\Lopto G_{ss}
+ G_{vu}\!\cdot\!\Pop G_{ss} \big) ,
\end{align*}
which effectively determines the characteristic numbers of curves in
$\P^1\times\P^1$ with one triple contact to a $(1,1)$-curve, from the usual
characteristic numbers.
\end{eks}

\small


\end{document}